


\magnification=1100
\overfullrule0pt

\input amssym.def


\def\qed{\hbox{\hskip 1pt\vrule width4pt height 6pt depth1.5pt \hskip 1pt}}
\def\mapright#1{\smash{
   \mathop{\longrightarrow}\limits^{#1}}}

\def\mapdown#1{\Big\downarrow
   \rlap{$\vcenter{\hbox{$\scriptstyle#1$}}$}}

\def\CC{{\Bbb C}}
\def\PP{{\Bbb P}}
\def\QQ{{\Bbb Q}}
\def\RR{{\Bbb R}}
\def\ZZ{{\Bbb Z}}
\def\cI{{\cal I}}
\def\cL{{\cal L}}
\def\cO{{\cal O}}
\def\cT{{\cal T}}

\def\gg{{\goth g}}
\def\gh{{\goth h}}

\def\Lie{{\rm Lie}}
\def\pt{{\rm pt}}



\font\smallcaps=cmcsc10
\font\titlefont=cmr10 scaled \magstep1

\font\sectionfont=cmbx10
\font\tinyrm=cmr10 at 8pt


\newcount\sectno
\newcount\subsectno
\newcount\resultno

\def\section #1. #2\par{
\sectno=#1
\resultno=0
\bigskip\noindent{\sectionfont #1.  #2}~\medbreak}

\def\subsection #1\par{\bigskip\noindent{\it  #1} \medbreak}


\def\prop{ 
\bigskip\noindent{\bf Proposition. }\sl}

\def\remark{ \global\advance\resultno by 1
\bigskip\noindent{\bf Remark \the\sectno.\the\resultno. }}
\def\example{ \global\advance\resultno by 1
\bigskip\noindent{\bf Example \the\sectno.\the\resultno. }}
\def\cor{
\bigskip\noindent{\bf Corollary. }\sl}
\def\thm{
\bigskip\noindent{\bf Theorem. }\sl}
\def\defn{ \global\advance\resultno by 1
\bigskip\noindent{\it Definition \the\resultno. }\slrm}
\def\endthm{\rm\bigskip}

\def\pf{\rm\bigskip\noindent{\it Proof. }}
\def\pfsk{\rm\bigskip\noindent{\it Sketch of proof. }}
\def\endpf{\qed\hfil\bigskip}


\def\formula{\global\advance\resultno by 1
\eqno{(\the\resultno)}}
\def\formulano{\global\advance\resultno by 1 (\the\resultno)}
\def\tableno{\global\advance\resultno by 1
\the\resultno. }
\def\lformula{\global\advance\resultno by 1
\leqno(\the\resultno)}

\def\monthname {\ifcase\month\or January\or February\or March\or April\or
May\or June\or
July\or August\or September\or October\or November\or December\fi}

\newcount\mins  \newcount\hours  \hours=\time \mins=\time
\def\now{\divide\hours by60 \multiply\hours by60 \advance\mins by-\hours
     \divide\hours by60         
     \ifnum\hours>12 \advance\hours by-12
       \number\hours:\ifnum\mins<10 0\fi\number\mins\ P.M.\else
       \number\hours:\ifnum\mins<10 0\fi\number\mins\ A.M.\fi}


\nopagenumbers
\def\runningtitle{\smallcaps a pieri-chevalley formula}
\headline={\ifnum\pageno>1\eoheadline\else\firstheadline\fi}
\def\names{\smallcaps h. pittie\quad and\quad a. ram}
\def\firstheadline{\noindent Preliminary Draft \hfill  \today}
\def\firstheadline{}
\def\eoheadline{\ifodd\pageno\oddheadline\else\evenheadline\fi}
\def\oddheadline{\tenrm\hfil\runningtitle\hfil\folio}
\def\evenheadline{\tenrm \folio\hfil{\names}\hfil}

\vphantom{$ $}  
\vskip.75truein

\centerline{\titlefont A Pieri-Chevalley formula in the }
\smallskip\centerline{\titlefont K-theory of a G/B-bundle}
\bigskip
\centerline{\rm Harsh Pittie}
\centerline{Department of Mathematics, Graduate Center}
\centerline{City University of New York}
\centerline{New York, NY 10036}
\bigskip
\centerline{\rm Arun Ram${}^\ast$}
\centerline{Department of Mathematics}
\centerline{Princeton University}
\centerline{Princeton, NJ 08544}
\centerline{{\tt rama@math.princeton.edu}}

\footnote{}{\tinyrm ${}^\ast$ Research supported in part by National
Science Foundation grant DMS-9622985.}

\bigskip


Let $G$ be a complex, semisimple, simply-connected
algebraic group and $B \subseteq G$ a Borel subgroup.
We fix a smooth closed complex projective variety $X$ and a principal
algebraic $B$-bundle over it:
$B \longrightarrow E \mapright{\pi} X$.  For
any complex algebraic variety $F$ with a left algebraic $B$-action, we
denote by $E (F)$ the total space of the associated fibre bundle
with fibre $F$.  Thus $E(F) = E \times_B F$
and the projection to
$X$ is obtained from projection on the first factor.

Fix a maximal torus $T\subseteq B$ and let $W$ be its Weyl group.
For each $w\in W$ the Bruhat cell $Y_w^\circ=BwB\subseteq G/B$ and
the Schubert variety $Y_w=\overline{BwB}\subseteq G/B$ are $B$-stable subsets of
$G/B$ so we have inclusions of bundles $E(Y_w^\circ)\subseteq 
E(Y_w)\subseteq E(G/B)$.  The closed subvarieties $\Omega_w=E(Y_w)$
determine classes $[\cO_{\Omega_w}]$ in 
$K(E(G/B))$\footnote{${}^\dagger$}{For any smooth variety $V$, 
$K(V)$ is the Grothendieck ring of coherent $\cO_V$-modules.}.
In fact, by a well known result of Grothendieck, these classes
form a $K(X)$-basis for $K(E(G/B))$.  On $E(G/B))$ we also have ``homogeneous''
line bundles associated to irreducible representations of $B$ (see below).
The main result of this announcement is a formula for the tensor
product of the class of a homogeneous line bundle with a 
Schubert class, expressed as a $K(X)$-linear combination of Schubert
classes.  

We believe that this formula is the most general uniform result in the
intersection theory of Schubert classes: it is related to a recent result
of Fulton and Lascoux [FL] who presented a similar formula for a
$GL_n(\CC)/B$-bundle.  Indeed, in this case, their formula and ours
coincide once one knows how to translate between their 
combinatorics with tableaux and ours with Littelmann paths.  O. Mathieu has also
proved the positivity which is implied by our formula, see [FP, p. 101]. 
Applying the Chern character to our formula, and equating the lowest order terms
we obtain a relative version of the result of Chevalley [Ch] alluded to in the
title of this paper.

The ring $K(E(G/B))$ is a $K(X)$-module via the map 
$\pi^*\colon K(X)\to K(E(G/B))$.  Since $G/B$ has a unique fixed point 
for the $B$-action there is a canonical
section $\sigma\colon X\to E(G/B)$ of the bundle $E(G/B)$.
Consider the diagram
$$
\matrix{
G& \mapright{} &E(G)& \mapright{} &X \cr
\mapdown{}&& \mapdown{\rho}&& \Vert \cr
G/B & \mapright{} &E(G/B)& \mapright{} &X \cr}
$$
where the vertical maps are quotients by the right action of
$B$ on $G$; precisely, $E(G/B) \simeq (E\times _B G) / B$.  Thus
$\rho$ is the projection map of a principal $B$-bundle over 
$E(G/B)$.

There are two vector bundles naturally associated to each $B$-module
$V$:
$$
E(V) \longrightarrow X, \qquad\hbox{and}\qquad
E_G(V)=E(G)\times_B V \longrightarrow E(G/B),$$
where the projection map for the latter of these is via $\rho$.
This assignment $V \mapsto E_G (V)$ of $B$-modules to vector
bundles over $E(G/B)$ preserves direct sums and tensor products,
and hence induces a ring homomorphism 
$R(B) \mapright{\phi} K (E(G/B))$,
where $R(B)$ is the representation ring of $B$. 
By construction $\sigma^*(E_G(V))=E(V)$ as vector bundles on $X$.
One also checks that if $V$ is the restriction
of a $G$ module, then $E_G(V)=\pi^*(\sigma^*(E_G(V)))$. 
Thus we have a commutative diagram
$$\matrix{
R(G) &-\ -\ \to &K(X) \cr
\mapdown{{\rm res}} &&\mapdown{\pi^*} \cr
R(B) &\longrightarrow &K(E(G/B))\cr
}$$
and a map
$$K(X)\otimes_{R(G)} R(B)\quad \mapright{\pi^*\otimes\phi}\quad K(E(G/B)),$$
where $R(G)$ is the representation ring of $G$ and the $R(G)$-action
on $K(X)$ is given by the map $V\mapsto E(V)$.

Let $P$ be the weight lattice of 
$\gg=\Lie(G)$. Then $R(B)=R(T)\cong \ZZ[P]$, the group algebra of $P$,
and $R(G)=R(T)^W$.  If $\lambda\in P$, let $e^\lambda$ be the corresponding
element of $R(T)$ and define
$$x^\lambda = E(e^\lambda)\in K(X)
\qquad\hbox{and}\qquad
y^\lambda = E_G(e^\lambda)\in K(E(G/B)).\formula$$
The statement that
$E_G(V)=\pi^*(\sigma^*(E_G(V)))$ if $V$ is a $G$-module is equivalent
to the statement that, in $K(E(G/B))$, 
$$\chi(x) = E(\chi)
\quad\hbox{is equal to}\quad
\chi(y)=E_G(\chi),
\quad\hbox{for all $\chi\in R(T)^W$.}$$
We recall from [P] that $R(T)$ is a free $R(G)$-module of rank
$|W|$, and that $R(T)\otimes _{R(G)} \Bbb Z\longrightarrow K(G/B)$
is an isomorphism.\footnote{${}^\dagger$}{The discussion in [P] is entirely in
terms of compact groups and the $K$-theory of $C^\infty$ vector
bundles;  with trivial modifications the results hold
in the present context also.} 
According to Steinberg, [S] there is an $R(G)$ basis of
$R(T)$ of the form $\{e^{\varepsilon_w}\ |\ w\in W\}$, where the $\varepsilon_w$
are certain specific elements of $P$.
Since the set $\{y^{\varepsilon_w}\ |\ w\in W\}$
is a set of globally defined elements in $K(E(G/B))$ which behaves 
properly under restriction, and which forms a basis locally,
it follows from standard yoga that it is also 
is a $K(X)$-basis for $K(E(G/B))$.  Thus the map
$K(X)\otimes_{R(G)} R(T) \longrightarrow K(E(G/B))$ is an isomorphism and
$$K(E(G/B))\cong {K(X)\otimes R(T)\over \cI},\formula$$
where $\cI$ is the ideal in $K(X)\otimes R(T)$ generated by the
set $\{ \chi(x)\otimes 1 -1\otimes \chi\ |\ \chi\in R(T)^W\}$.

Define a $W$-action on $K(X)\otimes R(T)$ as the $K(X)$-linear extension of the
action given by
$$wy^\lambda= y^{w\lambda},\qquad\hbox{for $w\in W$, $\lambda\in P$.}$$
This action descends to an action on $K(E(G/B)$ since the generators of 
the ideal $\cI$ are $W$-invariants for this action. 
Using this $W$-action on $K(E(G/B))$, we can
define the analogues of BGG-operators in this context.
Such operators were defined in the ``absolute case''
($X=\pt$) by Demazure, in $K_T (G/B)$ by Kostant and Kumar
[KK], and finally by Fulton and Lascoux [FL] when $G=  SL(n, \Bbb C)$.
To make the definition, let $\alpha$ be a positive root with respect
to the pair $(B,T)$ and let $s_\alpha \in W$ be the corresponding
reflection.
Define $T_\alpha\colon R(T)\to R(T)$ by defining
$$T_\alpha (e^\lambda ) 
= (e^{\lambda + \alpha} - s_\alpha (e^\lambda) ) / (e^\alpha -1)$$ 
and extending $\ZZ$-linearly.  Since $T_\alpha$ fixes elements of
$R(T)^W$ this operation can be extended $K(X)$-linearly to a well defined
operator on $K(E(G/B))$.

Now fix a simple system of roots $\alpha_1 , \ldots , \alpha _\ell$
for $(B,T)$ and let $P_j$ be the minimal parabolic subgroup 
corresponding to $\alpha_j$; this is the closed connected
subgroup of $G$ whose Lie algebra $\frak{p}_j$ is spanned by the
Lie algebra $\frak b$ of $B$ and the root space $\gg_{-\alpha_j}$.
Let $f_j : E(G/B) \longrightarrow E (G/P_j)$ be the 
projection induced from the $B$-equivariant ${\Bbb P}^1$-bundle
$G/B \longrightarrow G/P_j$ (the canonical projection). 
The following result explains the geometric significance of the
operators $T_{\alpha_j}$ (henceforth abbreviated as $T_j$).
P. Deligne pointed out an error in the proof of (a) below in an earlier
version of this preprint.  We are grateful to him for pointing this
out and have corrected the argument.

\prop
With the notation as above
\smallskip
\item{(a)} $\displaystyle{
(f_j)^! \circ (f_j)_! ([\cO_{\Omega_w}]) =
\cases{
[\cO_{\Omega_{ w s_j}} ], &if $\ell (w s_j ) > \ell (w)$, \cr
\cr
[\cO_{\Omega_w}], &if $\ell(w s_j) < \ell (w)$. \cr}
}
$
\smallskip
\item{(b)} \ For any element $x \in K(E(G/B))$,
$(f_j)^! \circ (f_j)_!(x) = T_j (x)$.
\pf
(a) Let $\bar w = \{w,ws_j\}$ be the coset of $w$ relative to 
$\langle s_j\rangle$.  The essential point is to prove the following
two equations
$$f_! [\cO_{\Omega_v}] = [\cO_{\Omega_{\bar w}}],
\qquad\hbox{for $v\in \bar w$,}$$
where $\Omega_{\bar w}\subseteq E(G/P_j)$ is the relative Schubert variety
constructed from $Y_{\bar w}\subseteq E(G/P_j)$.  In turn, these equations
will follow from the isomorphisms
$$\hbox{(i)}\quad f_*(\cO_{\Omega_v})=\cO_{\Omega_{\bar w}},
\qquad\quad
\hbox{(ii)}\quad R^qf_*(\cO_{\Omega_v}),
\quad\hbox{for $q>0$, $v\in \bar w$.}$$

To prove (i) and (ii) relabel the elements of $\bar w$ as $w'$
and $w''$ where $w'<w''$.  Then $f\colon \Omega_{w'}\to \Omega_{\bar w}$
is birational, and since the varieties in question have at worst
rational singularities, (i) and (ii) for $w'$ follow from known
arguments.  Secondly $f:\Omega_{w''}\to \Omega_{\bar w}$ is a
$\PP^1$-bundle, so (i) and (ii) are standard.  Finally
$f^![\cO_{\Omega_{\bar w}}] = [\cO_{\Omega_{w''}}]$ follows because,
in this case, $f$ is the projection of a $\PP^1$-bundle.
\smallskip\noindent
(b) There is a 2-dimensional algebraic
vector bundle $E_j \longrightarrow E(G/P_j)$ associated to a
2-dimensional representation of $P_j$.  Its projectivization
is $E(G/B)$, i.e.
$\Bbb P (E_j) \simeq E(G/B)$ as bundles over
$E(G/P_j)$. It follows that $K(E(G/B))$ is a free module over
$K(E(G/P_j))$ on two generators, 1 and $L_{\omega_j}$, where $\omega_j$ is the
$j^{\rm th}$ fundamental weight.  Since both sides are $K(X)$-linear,
it suffices to check the assertion for $1$ and $L_{\omega_j}$, and
this reduces to the ``same'' computation as in the absolute case.
\endpf

The operators $T_i$, $1\le i\le \ell$, satisfy $T_i^2=T_i$ and the generalized
braid relations. For each $w\in W$, let $w=s_{i_1}\cdots
s_{i_p}$ be a reduced word for $w$ and define $T_w=T_{i_1}\cdots T_{i_p}$.
Since the $T_i$ satisfy the braid relations, the operators $T_w$
are well defined and, by the above Proposition,
$$[\cO_{\Omega_w}]=T_{w^{-1}}[\cO_{\Omega_1}],
\qquad\hbox{for $w\in W$}.\formula$$

For each $\lambda\in P$ let $Y^\lambda$ be
the ``left multiplication'' operator on $K(E(G/B))$ defined by
$Y^\lambda(x) = y^\lambda x$.  Since $[\cO_{\Omega_1}]=\sigma(X)$,
$$Y^\lambda [\cO_{\Omega_1}] = x^\lambda [\cO_{\Omega_1}],\formula$$
where $x^\lambda$ is as in (1).  As operators on $K(E(G/B))$,
$$
Y^\lambda T_i = 
T_i Y^{s_i\lambda} + 
{Y^\lambda-Y^{s_i\lambda}\over 1-Y^{-\alpha_i}}, \formula$$
where the second term is always viewed as a linear combination 
$Y^\mu$, $\mu\in P$.
We will iterate this formula to obtain an expansion of the product
$e^\lambda[\cO_{X_w}]$ in $K(G/B)$ in terms of the $K(X)$-basis 
$\{[\cO_{X_v}]\ |\ v\in W\}$ of $K(E(G/B)$.
The path model of P. Littelmann [Li] is exactly what is needed for 
controlling the resulting expansion.

Let $\gh^*=\RR\otimes P$ be the real span of the weight lattice.
A path in $\gh^*$ is a piecewise linear map 
$\pi\colon [0,1]\to \gh^*$ such that $\pi(0)=0$.  
P. Littelmann [Li] defined {\it root operators} $f_1,\ldots, f_\ell$ which act 
on the paths. The action of a root operator $f_i$ on a path $\pi$ either 
produces another path or returns $0$. 

Let $\lambda$ be a dominant integral weight and let
$W_\lambda$ be the stabilizer of $\lambda$.  The cosets in $W/W_\lambda$ are
partially ordered by the Bruhat-Chevalley order.
Let $\pi_\lambda$ be the path
given by
$$\pi_\lambda(t)=t\lambda,\enspace 0\le t\le 1,
\quad\hbox{and let}\quad
\cT^\lambda = \{ f_{i_1}f_{i_2}\cdots f_{i_l}\pi_\lambda \}$$ 
be the set of all paths obtained by applying sequences of
root operators $f_i=f_{\alpha_i}$, $1\le i\le \ell$ to $\pi_\lambda$.
Each path $\pi\in \cT^\lambda$
can be encoded with a pair of sequences
$$\matrix{
\hfill\vec \tau&=(\tau_1>\tau_2>\cdots >\tau_r),\hfill
&\tau_i\in W/W_\lambda, \qquad \hbox{and}\hfill\cr
\hfill\vec a&=(0=a_0<a_1<a_2<\cdots<a_r=1),\hfill
&\qquad a_i\in \QQ,\hfill\cr
}$$
so that $\pi$ is given by
$$\pi(t)=(t-a_{j-1})\tau_j\lambda+\sum_{i=1}^{j-1} (a_i-a_{i-1})\tau_i\lambda,
\qquad\hbox{for $a_{j-1}\le t\le a_j$.}
$$
The {\it initial direction} of $\pi$ is $\iota(\pi)=\tau_1$ and the {\it endpoint}
of $\pi$ is $\pi(1)\in \gh^*$.

Fix $w\in W$, let $\bar w=wW_\lambda\in W/W_\lambda$ and assume that 
$\pi$ is a path in the set
$$\cT^\lambda_{\le w} = \{ \pi\in \cT^\lambda \ |\ \iota(\pi)\le \bar w\}.$$
A {\it maximal lift of $\vec\tau$ with respect to $w$} is a choice of 
representatives $t_i\in W$ of the cosets $\tau_i$ such that 
$w\ge t_1>\cdots>t_r$ and each $t_i$ is maximal in Bruhat order such that
$t_{i-1}>t_i$.   The {\it final direction} of $\pi$ with respect to $w$ is
$$v(\pi,w)=t_r,$$
where
$w\ge t_1>\cdots>t_r$ is a maximal lift of $\tau_1>\ldots>\tau_r$
with respect to $w$.

\thm  Let $\lambda$ be a dominant integral weight and let $w\in W$.  
Then
$$Y^\lambda T_{w^{-1}} = 
\sum_{\eta\in \cT^\lambda_{\le w}} T_{v(\eta,w)^{-1}} Y^{\eta(1)}
$$
as operators on $K(E(G/B))$.
\pfsk
Fix a simple root $\alpha_i$.  Every path is in
a unique {\it $\alpha_i$-string} of paths
$$S_{\alpha_i}(\pi)=\{f_i^m\pi,\ldots, f_i^2\pi,f_i\pi,\pi\},$$
where $f_i^m\pi=0$ and there does not exist any path $\eta$ such that
$f_i\eta=\pi$.  In a manner similar to that of [Li, Lemma 5.3] one easily
shows that, for any $\alpha_i$-string $S_{\alpha_i}(\pi)$,
$$\sum_{\eta\in S_{\alpha_i}(\pi)} T_{v(\eta,w)^{-1}}Y^{\eta(1)}
=T_{v(\pi,w^{-1})}Y^{\pi(1)}T_i.$$
Given these facts, the proof of the Theorem follows the same lines as the proof
of the Demazure character formula given in [Li, 5.5]. 
\endpf

By applying the formula in the Theorem to the element $[\cO_{\Omega_1}]\in
K(E(G/B))$ and using (3) and (4) we obtain the following.

\cor  Let $\lambda$ be a dominant integral weight
and let $w\in W$.  In $K(E(G/B))$,
$$y^\lambda [\cO_{\Omega_w}] = 
\sum_{\eta\in \cT^\lambda_w} [\cO_{\Omega_{v(\eta,w)}}]x^{\eta(1)}.
$$
\endthm

\bigskip\noindent
{\bf Acknowledgements.}
We are grateful to many people for comments,
suggestions and encouragement: we would particularly like to thank
Jim Carlson, Mark Green, Shrawan Kumar, Bob MacPherson, Ted Shifrin 
and Al Vasquez.

\bigskip\bigskip
\centerline{\smallcaps References}

\bigskip
\item{[Ch]} {\smallcaps C.\ Chevalley}, {\it Sur les decompositions cellulaires
des espaces $G/B$},  in {\sl Algebraic Groups and their Generalizations:
Classical Methods}, W. Haboush and B. Parshall eds.,
Proc. Symp. Pure Math., Vol. {\bf 56} Pt. 1, Amer. Math. Soc. (1994), 1--23.

\medskip
\item{[FL]} {\smallcaps W.\ Fulton and A.\ Lascoux},
{\it A Pieri formula in the Grothendieck ring of a flag bundle},
Duke Math. J. {\bf 76} (1994), 711--729.

\medskip
\item{[FL]} {\smallcaps W.\ Fulton and P. Pragacz},
{\sl Schubert varieties and degeneracy loci},
Lecture Notes in Math. {\bf 1689}, Springer-Verlag, Berlin 1998.

\medskip
\item{[KK]} {\smallcaps B.\ Kostant and S.\ Kumar}, 
{\it $T$-equivariant K-theory of generalized flag varieties}, J. Differential
Geom. {\bf 32} (1990), 549--603. 

\medskip
\item{[Li]} {\smallcaps P.\ Littelmann},
{\it  A Littlewood-Richardson rule for symmetrizable Kac-Moody algebras},
Invent. Math. {\bf 116} (1994), 329-346.

\medskip
\item{[P]} {\smallcaps H.\ Pittie}, 
{\it Homogeneous vector bundles over homogeneous spaces},
Topology {\bf 11} (1972), 199--203.

\medskip
\item{[S]} {\smallcaps R.\ Steinberg}, 
{\it On a theorem of Pittie},
Topology {\bf 14} (1975), 173-177.

\vfill\eject
\end

In the total space $E(G/B)$ we can define Bruhat ``cells''
and Schubert varieties ``parametrized by $X$'' as follows.
Fix a maximal torus $T \subseteq B$ and let $W$ be its Weyl group.
Note that since the identity coset $1 \cdot B \in G/B$
is a fixed point of the $B$-action, we obtain an
algebraic section $\sigma: X \longrightarrow E (G/B)$.
By definition, $\sigma (X) = \Omega _1$ is the Schubert variety
corresponding to $1 \in W$.  For general $w \in W$ let
$U \subseteq X$ be an open set over which
$\pi : E \longrightarrow X$, (and hence also
the associated bundle $E(G/B) \longrightarrow X)$ is trivial.
Then we define the $B$-orbit of $w$ along $U$ as the set of 
pairs $\{ (u, b w \sigma (u)) : u \in U , b \in B \}$
as a subset of $U \times G/B$.  Since the transition functions of
$E(G/B) \longrightarrow X$ lie in $B$, these local  $B$-orbits
fit together to give a variety $\Omega _w^\circ \subseteq E (G/B)$,
which is the Bruhat cell corresponding to $w \in W$.
Note that for any $x \in X$, $\pi^{-1} (x) \cap \Omega_ w^\circ$
is the usual Bruhat cell in $\pi^{-1} (x) \simeq G/B$:
in particular, $dim_{\Bbb C} \Omega_w^\circ = \ell (w) + dim_{\Bbb C}X$,
where $\ell(w) =$ length of $w$.  Then the corresponding
Schubert cycle $\Omega_w =$ closure of $\Omega_w^\circ$ in $E(G/B)$.

Let $\cO_{\Omega_w}$ be the structure sheaf of $\Omega_w$: we use the
same notation for the corresponding sheaf on $E(G/B)$, which is
$\cO_{\Omega_w}$ extended by zero on $E(G/B) - \Omega_w$.
A well known result of Grothendieck [C] implies that
$K(E(G/B))$ is a free $K(X)$-module (under $\pi^!)$ with
basis $\{ [\cO_{\Omega w}]: w \in W \}$.  Here and henceforth,
$K(Y)$ denotes the  Grothendieck ring of coherent $\cO_Y$-modules.
(Most of the $Y$ we will have occasion to use will be smooth;
hence $K(Y)$ is canonically isomorphic to the Grothendieck ring of
{\it locally free} $\cO_Y$-modules.)

Our main formula gives the tensor product of a (positive) element of
this basis with $[\cO_{\Omega _w }]$, expressed as a $K(X)$-linear
combination of $\{[\cO_{\Omega_w}]: w \in W\}$.  To describe this
construction,
take note of the diagram

The principal bundle $\rho $ and hence also $E_G (V)$ is a bundle
``along the fibres'' of $\pi : E(G/B) \longrightarrow X$,
in the sense that locally the transition functions are
constant in the horizontal direction, if $U \subseteq X$
is an open set over which
$\pi^{-1} (U) \simeq U \times B$, so that
$E(G/B) |_U \simeq U \times G/B$, then over
$\Phi = E (G/B)|_U$,
$$
E_G (V) | _\Phi = U \times (G \times _BV)
$$
with the transition functions being $1_U$ on the first factor.

is a trivial bundle over
$E(G/B)$ of rank equal to $\dim_{\Bbb C}V$.  Hence the
map $\varphi$ factors to give 
$\overline{\varphi} : R(B) \otimes _{R(G)} \Bbb Z \longrightarrow K(E(G/B))$.
Using Lie's theorem $R(B)$ is canonically isomorphic to $R(T)$,
and so we can rewrite $\overline{\varphi}$ as
$R(T) \otimes _{R(G)} \Bbb Z \longrightarrow K(E(G/B))$.
This notation has the advantage that it will enable us to define
(formally) an action of $W$ on $K(E(G/B))$.

in fact from the discussion above this isomorphism is 
$\overline{\varphi}$ composed with restriction to a fibre.
Thus the restriction map 
$K(E(G/B)) \mapright{i^*} K(G/B)$ admits
an additive retract.  The importance of this observation will become
apparent presently.

:  that is, if $w \in W$ is expressed
as a matrix in terms of a basis $\{\chi_1 , \ldots , \chi_n \}$
it lies in $GL(n, \Bbb Z)$
 (rather than $K(X)$)!

The $R(G)$-action of $W$ on $R(T)$ gives an action of $W$ on 
$R(T) \otimes  _{R(G)} \Bbb Z$
which is $\Bbb Z$ linear.  We use $\overline{\varphi}$ to
transport this to a $GL(n, \Bbb Z)$-action on the $\Bbb Z$-span
of $\cL$: and then, since $\cL$ is a 
$K(X)$ basis for $K(E(G/B))$ we obtain a
$K(X)$-linear action of $W$
on $K(E(G/B))$, whose matrix entries in the basis $\cL$ are in
$\Bbb Z$.